\documentclass[12pt]{article}
\usepackage{amsmath,amssymb}

\newcommand{\Var}{\ensuremath{\mathcal{V}_{\mathbb{C}}}}

\def\vv{{\underline{v}}}
\def\tt{{\underline{t}}}

\def\MM{\underline{M}}

\def\NN{\underline{N}}
\def\kk{{\underline{k}}}
\def\nuunder{\underline{\nu}}
\def\1{\underline{1}}

\def\P{\mathbb P}
\def\N{\mathbb N}
\def\LLL{\mathbb L}
\def\Z{\mathbb Z}

\def\C{\mathbb C}
\def\Ex{\mathbf{Exp}}

\def\S{{\mathcal S}}
\def\RR{{\mathcal R}}
\def\OO{{\mathcal O}}
\def\X{{\mathcal X}}
\def\D{{\mathcal D}}
\def\B{{\mathcal B}}
\def\A{{\mathcal A}}
\def\LL{{\mathcal L}}
\def\m{{\frak m}}

\newtheorem{theorem}{Theorem}

\newtheorem{proposition}{Proposition}

\newenvironment{definition}
{\smallskip\noindent{\bf Definition\/}:}{\smallskip\par}
\newenvironment{corollary}
{\smallskip\noindent{\bf Corollary\/}.}{\smallskip\par}

\newenvironment{proof}
{\noindent{\bf Proof\/}.}{{ $\square$}\smallskip\par}

\title{Integration over spaces of non-parametrized arcs
and motivic versions of the monodromy zeta function
\footnote{Math. Subject Class. 14B05, 32S05, 58K10}
}
\author{S.M.~Gusein-Zade \thanks{Partially supported by the grants
RFBR--04--01--00762,
NSh--1972.2003.1.
Address: Moscow State University,
Faculty of Mathematics and Mechanics, Moscow, 119992, Russia.
E-mail: sabir\symbol{'100}mccme.ru} \and
I.~Luengo
\and A.~Melle--Hern\'andez \thanks{The last two authors were partially
supported by the grant BFM2001--1488--C02--01. Address:
University Complutense de Madrid, Dept. of Algebra,
Madrid, 28040, Spain.
E-mail: iluengo\symbol{'100}mat.ucm.es, amelle\symbol{'100}mat.ucm.es}}
\date{}
\begin{document}
\def\eps{\varepsilon}

\maketitle

\begin{abstract}
We elaborate notions of integration over the space of arcs factorized by the natural $\C^*$-action and over the space of non-parametrized arcs (branches). There are offered two motivic versions of the zeta function of the classical monodromy transformation of a germ of an analytic function on $\C^d$. We indicate a direct formula which connects the naive motivic zeta function of J.~Denef and F.~Loeser with the classical monodromy zeta function.
\end{abstract}
\medskip

\section*{Introduction.}
The notion of motivic integration invented by M.~Kontsevich and developed by V.~Batyrev, J.~Denef, F.~Loeser et al. (see, e.g., \cite{DL2, DL3, L}) is an
analogue of $p$-adic integration. It also can be considered as a generalization of the notion of integration with respect to the Euler characteristic (see \cite{Viro}) in two directions. First, instead of the usual Euler characteristic (with values in the ring of integers $\Z$) one considers the generalized (universal) Euler characteristic with values in the Grothendieck ring $K_0(\Var)$ of complex algebraic varieties or/and in a modification (localization, completion) of it. Second, instead of integration over, say, a (finite dimensional) algebraic variety one integrates over the infinite dimensional space of arcs. This notion, in particular, gives a possibility to construct (or to define) motivic versions of some classical invariants of varieties or of singularities. The notion of a motivic version of an invariant is not well defined. There are only two obvious requirements: such a version should be an invariant itself and it should specialize to the classical one under a corresponding additive invariant: Euler characteristic, Hodge--Deligne polynomial, ... Sometimes such invariants can be defined as certain integrals with respect to the universal Euler characteristic. However, one can meet the following problem.

On the space of arcs, there is a natural $\C^*$-action defined by $a*\varphi(\tau)=\varphi(a\tau)$. Majority of natural constructible functions on the space of arcs which could participate in a definition of an invariant (say, the order of a fixed function along an arc) are invariant with respect to this action. The integral of such a function over the space of arcs with respect to the universal Euler characteristic is divisible by the class $\LLL-1$ of the punctured complex line. Therefore the specialization of this integral by the usual Euler characteristic morphism is equal to zero and no motivic version of an usual invariant can be constructed this way (as such an integral). E.g., the Euler characteristic of the "naive" zeta function of J.~Denef and F.~Loeser (see \cite{DL3}) is equal to zero.

To define motivic versions of (integer valued) invariants of singularities, one has "to kill" this $\C^*$-action. One can imagine several ways to do this. One way is to consider not the whole space of arcs but a certain subspace of it. This was made in \cite{DL4} where instead of the space of arcs 
%% there
it was considered its subspace consisting of arcs with the prescribed first coefficient of the Taylor expantion of the function under consideration along an arc.

As another possibility one can imagine substitution of integration over the space of arcs by integration over another (infinite-dimensional) space. For example, integration over the space of functions or over its projectivization appeared to be useful for some problems: \cite{CDG1}, \cite{CDG2}, {\dots} For instance one can factorize the space of arcs by this $\C^*$-action. One can say that this way is inspired by the notion of projectivization. (It is not really the projectivization since this $\C^*$-action is not free.) For integration over the space of functions it was used, e.g., in \cite{CDG2}.

It is also possible instead of arcs to try to consider branches, under what we understand arcs without parametrization. In other words we consider the space of arcs (space of maps $\varphi:(\C,0)\to W$) factorized by the group $Aut_{\C,0}$ of changes of the coordinate $(\C,0)\to(\C,0)$. One can say that this notion has more geometric meaning than arcs modulo the described $\C^*$-action, in particular, since this action itself depends on the choice of the coordinate in the source $(\C,0)$. Moreover it seems that integration over the set of branches should be more generalizable to possible constructions of integration over sets
of higher dimensional subspaces.

We shall be interested in constructing 
%% (defining)
a motivic version of the classical monodromy zeta function
$$
\zeta_f(t)=\prod\limits_{q\geq 0} \{ \det\,[\,id-t\,h_{*}
\raisebox{-0.5ex}{$\vert$}{}_{{H_q(V_f;\C)}}]\}^{(-1)^{q+1}}
$$
of the germ of a function $f:(\C^d,0)\to(\C,0)$ ($h:V_f\to V_f$ is the classical monodromy transformation of the germ $f$). By the A'Campo formula (\cite{A'C}) the zeta function $\zeta_f(t)$ can be written as the integral of the expression $(1-t^m)$ over the exceptional divisor of a resolution of the germ $f$ (in the group $1+t\cdot\Z[[t]]$ with respect to multiplication). Arcs on $\C^d$ correspond to points of the exceptional divisor. Thus one can substitute integration over the exceptional divisor by integration over the space of arcs. However, because of the described reason, the corresponding integral over the space of arcs degenerates to $1$ under the specialization to the usual Euler characteristic. To avoid this problem, we shall elaborate notions of integration over the space of arcs factorized by the $\C^*$-action and over the space of branches. Here we discuss these notions only in the smooth case. An integral over the space of arcs factorized by $\C^*$ can be considered as a well-defined division by $(\LLL-1)$ of the corresponding integral over the space of arcs itself (otherwise defined only up to torsion: see, e.g., \cite{DL1}). 

Another problem which can be met on this way is to give meaning to an expression of the form $(1-t^m)^{-\chi(Z)}$ in the A'Campo formula when the usual Euler characteristic $\chi(Z)$ is substituted by the universal one (i.e. to give sense to the expression $(1-t^m)^{-[Z]}$ for $[Z]$ from the Grothendieck ring $K_0(\Var)$ or from its localization $K_0(\Var)[\LLL^{-1}]$ by the class $\LLL$ of the complex affine line). A way to do this was the main result of \cite{MRL}.

The zeta function $\zeta_f(t)$ of the monodromy transformation can be obtained from the Denef--Loeser motivic Milnor fibre (see \cite{DL4}). The last one is the limit at infinity of the Igusa motivic zeta function. The motivic Milnor fibre can be described in terms of ramified coverings of components of the exceptional divisor of a resolution of the germ $f$. The motivic invariants which specialize to the zeta function $\zeta_f(t)$ offered here are expressed in terms of components of the exceptional divisor themselves, not in terms of their coverings.

The Grothendieck semiring $S_0(\Var)$ of complex quasi-projective varieties is the semigroup generated by isomorphism classes $[X]$ of such varieties modulo the relation $[X]=[X-Y]+[Y]$ for a Zariski closed subvariety $Y\subset X$; the multiplication is defined by Cartesian product: $[X_1]\cdot [X_2]=[X_1\times X_2]$. The Grothendieck ring $K_0(\Var)$ is the group generated by these classes with the same relation and the same multiplication. Let $\LLL
\in K_0(\Var)$ be the class of the complex affine line, and let $K_0(\Var)[\LLL^{-1}]$ be the localization of Grothendieck ring $K_0(\Var)$ with respect to $\LLL$. The class $[X]\in S_0(\Var)$ can be defined for any constructible subset $X$ (as $\sum [X_i]$ for a partition $X=\cup X_i$ of the set $X$ into a finite union of quasi-projective varieties). 

In what follows $\RR$ will denote one of the discussed (semi)rings. There is a natural (semi)ring homomorphism $\chi:\RR\to\Z$ which sends the class $[X]$ of a variety $X$ to the Euler characteristic $\chi(X)$ of the set $X$. For a series $A(t)=\sum\limits_{i=0}^\infty A_i\, t^i$ with coefficients $A_i$ from $\RR$ its specialization $\chi(A(t))$ under the Euler characteristic homomorphism is the power series $\sum\limits_{i=0}^\infty \chi(A_i)\,t^i\in\Z[[t]]$. 

Definition of the space of arcs on an algebraic variety and of the motivic measure on it can be found, e.g., in \cite{DL2, L} (here we use them only for smooth varieties).

\section{The space of arcs factorized by the $\C^*$-action.}
Let $\LL_0$ be the space of arcs on the affine space $\C^d$ at the origin, i.e. maps $\varphi:(\C,0)\to(\C^d,0)$ (in particular $\LL_0\subset\m\OO_{\C^d,0}^d$, where $\m$ is the maximal ideal in the ring $\OO_{\C^d,0}$ of germs of functions). For $n\ge 0$, let $\LL_0^n$ be the space of $n$-jets of arcs, i.e. the space $\LL_0/\m^n\LL_0$ of arcs truncated at the level $n$. There is a natural $\C^*$-action on the spaces $\LL_0$ and $\LL_0^n$ defined by $a*\varphi(\tau)=\varphi(a\tau)$ ($a\in\C^*=\C\setminus\{0\}$). Let $\LL_0^*:=\LL_0\setminus\{0\}$, $\LL_0^{*n}:=\LL_0^n\setminus\{0\}$ and let $\LL_0^*/\C^*$ and $\LL_0^{n*}/\C^*$ be the corresponding spaces factorized by the $\C^*$-action. The space $\LL_0^{n*}/\C^*$ is a (finite dimensional) projective variety. Therefore, for a constructible subset $Y$ of it, there is defined its generalized (universal) Euler characteristic $\chi_g(Y)=[Y]\in K_0(\Var)$. For $n\ge 0$, there exists a natural map $\pi_n:\LL_0^*/\C^*\to(\LL_0^{n*}/\C^*)\cup \{0\}$ and, for $n\ge m$, there exists a natural map $\pi_{n,m}:(\LL_0^{n*}/\C^*)\cup \{0\}\to(\LL_0^{m*}/\C^*)\cup \{0\}$. The last map is constructible. The space $\LL_0^{m*}/\C^*$ can be decomposed into finitely many constuctible subsets so that over each of them the map $\pi_{n,m}$ is a locally trivial fibration (in the Zariski topology) whose fibre is a vector space of dimension $d(n-m)$ factorized by a finite cyclic group action (the isotropy group of the corresponding jet). Since the class in $K_0(\Var)$ of a vector space of dimension $d$ factorized by a representation of a finite abelian group is equal to $\LLL^d$ (see \cite[Lemma 5.1]{L}) then, for a constructible subset $Y$ in $\LL_0^{m*}/\C^*$, one has $[\pi^{-1}_{n,m}(Y)]=\LLL^{d(n-m)}\cdot[Y]$. This inspires the following definitions.

\begin{definition} A subset $X\subset \LL_0^*/\C^*$ is called 
%% constructible 
{\em cylindric} if there exist $n\ge 0$ and a constructible subset $Y\subset\LL_0^{n*}/\C^*$ such that $X=\pi^{-1}_n(Y)$.
\end{definition}

\begin{definition} The {\em motivic measure} (or the {universal Euler characteristic}) of a constructible subset $X\subset \LL_0^*/\C^*$, $X=\pi^{-1}_n(Y)$ for $Y\subset\LL_0^{n*}/\C^*$, is $\chi_g(X):=[Y]\cdot\LLL^{-dn}\in K_0(\Var)[\LLL^{-1}]$.
\end{definition}

\begin{definition} A function $\psi:\LL_0^*/\C^*\to G$ with values in an abelian group $G$ is {\em constructible} if it has countably many values and, for each $a\in G$, $a\ne 0$, the level set $\psi^{-1}(a)$ is constructible.
\end{definition}

In the usual way (see, e.g., \cite{DL2, L}) one can define the integral $\int_{\LL_0^*/\C^*}\psi\,d\chi_g$ of the function $\psi$ with respect to the generalized Euler characteristic (the motivic measure), as $\sum\limits_{a\in G, a\ne0}\chi_g(\psi^{-1}(a))\cdot a$. (Warning: not all constructible functions are integrable since the sum of a series may have no sense in the group $G$.)

Let $p$ be the factorization map $\LL_0^*\to\LL_0^*/\C^*$. For a constructible subset $X\subset \LL_0^*/\C^*$, let $\widetilde X=p^{-1}(X)$ be the corresponding $\C^*$-invariant subset of the space $\LL_0$ of arcs. One can easily see that $\chi_g(\widetilde X)=(\LLL-1)\chi_g(X)$. This implies the following statement.

\begin{proposition}\label{prop1}
Let $\psi:\LL_0^*/\C^*\to G$ be a constructible integrable function and let $\widetilde\psi=\psi\circ p:\LL_0^*\to G$ be the corresponding $\C^*$-invariant function on the space of arcs. Then the function $\widetilde\psi$ is integrable and 
$$
\int_{\LL_0}\widetilde\psi\, d\chi_g=(\LLL-1)\int_{\LL_0^*/\C^*}\psi\, d\chi_g.
$$
\end{proposition}

This means that integrals over the space $\LL_0^*/\C^*$ can be considered as well defined versions of the corresponding integrals over the space of arcs itself divided by $(\LLL-1)$. It is not clear that such division is well defined in the ring $K_0(\Var)[\LLL^{-1}]$. Usually it can be made formally when such an integral is computed (say, in terms of a resolution), however, either the result should be considered modulo elements from the annulator of $(\LLL-1)$, or one should proof that the result does not depend on a resolution. For instance in \cite{DL4} it was shown that the motivic Milnor fibre (introduced in \cite{DL1} up to $(\LLL-1)$-torsion) is well defined. Therefore integration over the space $\LL_0^*/\C^*$ of arcs modulo $\C^*$ can be considered as a formalization of this procedure.

\section{The space of branches on $(\C^d,0)$.}

Now we adapt the construction described above for the space of arcs factorized by the group $Aut_{\C,0}$ of local changes of the coordinate in $(\C, 0)$. For an arc $\varphi:(\C, 0)\to(\C^d, 0)$ and $h\in Aut_{\C,0}$ ($h:(\C, 0)\to(\C, 0)$), let $h*\varphi(\tau):=\varphi(h^{-1}(\tau))$. This defines an action of the group $Aut_{\C,0}$ on the space of arcs. 

\begin{definition} An orbit of such an action is called a {\em branch} (on $(\C^d, 0)$).
\end{definition}

The group $Aut_{\C,0}$ acts also on the jet space $\LL_0^n$. Moreover such an action coincides with the induced action of the group $Aut_{\C,0}^n$ of $n$-jets of coordinate changes on $\LL_0^n$. Let $\B_0:=\LL_0/Aut_{\C,0}$ and $\B_0^n:=\LL_0^n/Aut_{\C,0}$ be the factor spaces of this action, let $\B_0^{n*}:=\B_0^n\setminus\{0\}=\LL_0^{n*}/Aut_{\C,0}$. Here all of them are considered simply as sets without an additional structure. 

The jet space  $\LL_0^n$ has a natural filtration defined by powers of the maximal ideal $\m$ of the ring $\OO_{\C^d,0}$ of germs of functions at the origin in $\C^d:$
$$
\{0\}\subset \m^{n-1} \LL_0^n \subset \m^{n-2} \LL_0^n\subset \ldots \subset \m^{1} \LL_0^n \subset \m^{0} \LL_0^n=\LL_0^n.
$$
This filtration respects the action of the group $Aut_{\C,0}.$ Let 
$$
\LL_0^{n*}=\bigcup_{i=0}^{n-1} \m^{i}\LL_0^n \setminus \m^{i+1}\LL_0^n
$$
be the corresponding decomposition of the punctured jet space $\LL_0^{n*}$ of $\LL_0^n$. For each $i$,  the isotropy groups of points for the action of the group $Aut_{\C,0}^n$ on the space $\m^{i}\LL_0^n \setminus \m^{i+1}\LL_0^n$ have the same dimension (equal to $n-i$). The factor space $(\m^{i}\LL_0^n \setminus \m^{i+1}\LL_0^n)/Aut_{\C,0}^n$ is a quasi-projective variety.

\begin{definition} A subset $Y\subset \B_0^{n*}$ will be called
{\em constructible} if for each $i=0,\ldots, n-1$ the set $Y_i=Y\cap \left((\m^{i}\LL_0^n \setminus \m^{i+1}\LL_0^n)/Aut_{\C,0}^n\right)$ is a constructible subset of the set $(\m^{i}\LL_0^n \setminus \m^{i+1}\LL_0^n)/Aut_{\C,0}^n$. The {\em generalized Euler characteristic} $\chi_g(Y)=[Y]$ of the set $Y$ is the sum $\sum\limits_{i=0}^{n} [Y_i]$ of classes of its parts $Y_i$.
\end{definition}

Let $\pi_n:\B_0\to\B_0^n$ and, for $n\ge m$, $\pi_{n,m}:\B_0^n\to\B_0^m$ be the natural maps.
For $n\geq m$, there exists a stratification of $\B_0^{m*}$ such that over each stratum the map $\pi_{n,m}:\B_0^n\to \B_0^m$ is a locally trivial fibration (in the Zariski topology), whose fibre is the factor of a $(d-1)(n-m)$-dimensional vector space by a finite cyclic group action. Therefore, for a constructible subset $Y$ in $\B_0^{m*}$, $n\ge m$, one has $[\pi_{n,m}^{-1}(Y)]=\LLL^{(d-1)(n-m)}[Y]$. This inspires the following definitions.

\begin{definition} A subset $X\subset \B_0$ is called {\em cylindric} if there exist $n\ge 0$ and a constructible subset $Y\subset\B_0^{n*}\subset\B_0^n$ such that $X=\pi^{-1}_n(Y)$.
\end{definition}

\begin{definition} The {\em motivic measure} (or the {\em universal Euler characteristic}) of a cylindric subset $X\subset \B_0$, $X=\pi^{-1}_n(Y)$ for $Y\subset\B_0^{n*}$, is $\chi_g(X)=[Y]\cdot\LLL^{-(d-1)n}\in K_0(\Var)[\LLL^{-1}]$.
\end{definition}

This measure induces the corresponding notion of integration over the space of branches.

Let $p_b$ be the factorization map $\LL_0 \to \LL_0/Aut_{\C,0}$, let $p^n_b$ the same map $\LL_0^n \to \LL_0^n/Aut_{\C,0}^n$: 
\[\begin{array}{ccccc}  & \LL_0 & \stackrel{\pi_n}{\longrightarrow} & \LL_0^n & \\ &  p_b \, \downarrow &&  \downarrow \, p_b^n &\\ & \B_0 & \stackrel{\pi_n}{\longrightarrow} & \B_0^n\end{array}\]
For a cylindric subset $X\subset \B_0$, $X=\pi_n^{-1}(Y)$, $Y\subset\B_0^{n*}$, let 
${\widetilde X}=p_b^{-1}(X)$ be the corresponding $Aut_{\C,0}$-invariant set of arcs. Let $Y_i=Y\cap \left((\m^{i}\LL_0^n \setminus \m^{i+1}\LL_0^n)/Aut_{\C,0} \right)$ and ${\widetilde Y_i}=(p^n_b)^{-1}(Y_i)$, for $i=0,1,\ldots,n-1$. 
One has $Y=\bigcup_{i=0}^n Y_i.$
Let $X_i:=(\pi_n)^{-1}(Y_i)$ and ${\widetilde X}_i=(\pi_n\circ p_b)^{-1}(Y_i)$, then
$X=\bigcup_{i=0}^{n-1} X_i$ and ${\widetilde X}=\bigcup_{i=0}^{n-1} {\widetilde X}_i$. 
%be the corresponding decompositions of the sets $X$ and ${\widetilde X}.$ 
One has $[{\widetilde Y}_i]=(\LLL-1)\LLL^{n-i-1} [Y_i].$ Therefore 
$$
\chi_g({\widetilde X}_i)=\LLL^{-dn}[{\widetilde Y}_i]=(\LLL-1)\,\LLL^{n-i-1-dn} [Y_i]=(\LLL-1)\,\LLL^{-i-1} \chi_g({X}_i).
$$

This implies the following statement. Let $ord $ be the order function on the space of arcs: $ord(\varphi)=i+1$ if $\varphi \in \m^i\LL_0\setminus \m^{i+1}\LL_0$, i.e. $ord(\varphi)=i+1$ for $\varphi\in {\widetilde X}_i$ ($ord$ is an $Aut_{\C,0}$-invariant function on the space of arcs). Let $T^{ord}$ be the corresponding map $\LL_0 \to \Z[[T]]$ (where $\Z[[T]]$ is considered as an abelian group with respect to summation). For a constructible function $\psi :\B_0\to G$ ($G$ is an abelian group), let ${\widetilde \psi}=\psi \circ p_b:\LL_0\to G$ be the corresponding $Aut_{\C,0}$-invariant function on the space of arcs and let ${\widetilde \psi}^*={\widetilde \psi}\otimes T^{ord}$ be the corresponding function with values in $G\otimes_{\Z} \Z[[T]]=G[[T]]$. 

\begin{proposition}\label{prop2}
If the function $\psi$ is a constructible integrable function then the function ${\widetilde \psi}^*$ is integrable and 
$$
\int_{\LL_0}\widetilde\psi^*\, d\chi_g\,
\raisebox{-0.5ex}{$\vert_{{T\mapsto\,\LLL}}$}=(\LLL-1)\int_{\B_0}\psi\, d\chi_g.
$$
\end{proposition}

Thus an integral over the space of branches of a $G$-valued function is defined by a certain integral over the space of arcs but of a $G[[T]]$-valued function.

\begin{proof} Let  $X_a:=\psi^{-1}(a).$ Then 
\begin{eqnarray*}
(\LLL-1)\int_{\B_0}\psi\, d\chi_g&=&(\LLL-1)\sum\limits_{a\in G, a\ne 0} \chi_g(X_a) \cdot a=\\
= (\LLL-1)\sum\limits_{a\in G, a\ne 0} \left(\sum\limits_{i=0}^{n(a)-1} \chi_g(X_{a,i}) \right)  \cdot a&=& 
\sum\limits_{a\in G, a\ne 0} \left(\sum\limits_{i=0}^{n(a)-1} \LLL^{i+1}\chi_g({\widetilde X}_{a,i}) \right) \cdot a = \\
=\sum\limits_{a\in G, a\ne 0} \left(\sum\limits_{i=0}^{n(a)-1} T^{ord} \chi_g({\widetilde X}_{a,i}) \right)\raisebox{-0.5ex}{$\vert_{{T\mapsto\,\LLL}}$} \, \cdot a &=& \int_{\LL_0}\widetilde\psi^*\, d\chi_g\,\raisebox{-0.5ex}{$\vert_{{T\mapsto\,\LLL}}$}
\end{eqnarray*}
(since, for $\varphi\in {\widetilde X}_{a,i}$, $ord(\varphi)=i+1$).
\end{proof}   

\section{Integrals in terms of a resolution.}
Let $f:(\C^d,0)\to(\C,0)$ be a germ of an analytic function, and let $\pi:(\X,\D)\to(\C,0)$ be a resolution of the germ $f$, i.e. a proper modification of $(\C^d,0)$ which is an isomorphism outside of the zero-level set $\{f=0\}$, $\X$ is smooth and the exceptional divisor $\D=\pi^{-1}(0)$ and the total transform ${\mathcal E}=(f\circ\pi)^{-1}(0)$ of the zero-level set of $f$ are normal crossing divisors on $\X$. For $\varphi\in\LL_0$, let $v_f(\varphi)=ord_f(\varphi)$ be the order of the function $f$ on the arc $\varphi$, $v_f=ord_f:\LL_0 \to \Z\cup\{\infty\}.$ The function $v_f$ is  $Aut_{\C,0}$-invariant (and therefore $\C^*$-invariant).

Let ${\mathcal E}=\cup_{i\in I_0} E_i$ be the decomposition of the total transform $\mathcal E$ into the union of irreducible components, let $I_0=I_0^{'}\cup I_0^{''}$
where, for  $i\in I_0^{'}$ (respec. for $i\in I_0^{''}$), $E_i\subset {\D}$ (respec. $E_i$ is a component of the strict transform of $\{f=0\}).$ For $i\in I_0,$
let $N_i$ be the multiplicity of the lifting $f\circ \pi$ of the function $f$ to the space $\X$ of the resolution along the component $E_i$, let $\nu_i-1$ be the multiplicity of the $d$-form $\pi^*\,dx$ along the corresponding $E_i$ ($dx=dx_1\wedge \ldots \wedge dx_d$ is the volume form on $\C^d$).   For $i\in I_0,$
let $\stackrel{\circ}{E_i}:=E_i \setminus \cup_{j\neq i} E_j$ be "the smooth part" of the component $E_i$; for $I\subset I_0$, $I\ne\emptyset$, let $E_I:= \cap_{i\in I} E_i$, $\stackrel{\circ}{E_I}:= E_I\setminus \cup_{j\in I_0\setminus I} E_j$.

The argument of \cite[Theorem 2.2.1]{DL1} imply the following statement.

\begin{proposition}\label{prop3}
$$
\int_{\LL_0^*/\C^*} t^{v_f}\,d\chi_g= \sum_{I\subset I_0, \emptyset\ne I\not\subset I_0^{''}}
(\LLL-1)^{|I|-1} [\stackrel{\circ}{E_I}] \prod_{i\in I} \frac{\LLL^{-\nu_i}\,t^{N_i}}{1-\LLL^{-\nu_i}\,t^{N_i}}.
$$
\end{proposition}

Now suppose that the resolution $\pi:(\X,\D)\to (\C^d,0)$ factorizes through the blowing-up $\pi_0:(\X_0, \C\P^{d-1})\to (\C^d,0)$ at the origin in $(\C^d,0)$, i.e. $\pi=\pi_0\circ \pi',\,\pi':(\X,\D) \to (\X_0,\C\P^{d-1})$. For $i\in I_0,$ let $M_i$ be the multiplicity of the component $E_i$ in the divisor 
$\pi'^*(\C\P^{d-1}).$ The same arguments with those used in Proposition 2 give the following statement. 

\begin{proposition}\label{prop4}
$$
\int_{\B_0} t^{v_f}\, d\chi_g= \sum_{I\subset I_0, \emptyset\ne I\not\subset I_0^{''}}
(\LLL-1)^{|I|-1} [\stackrel{\circ}{E_I}] \prod_{i\in I} \frac{\LLL^{-\nu_i-M_i}\,t^{N_i}}{1-\LLL^{-\nu_i-M_i}\,t^{N_i}}.
$$
\end{proposition}

\section{Power structure over the Grothendieck ring of varieties.}
In what follows we shall be interested in integrals with respect to the motivic measure of a function $(1-t^{v_f})^{-1}$ whose values are considered as elements of the group $1+t\cdot K_0(\Var)[\LLL^{-1}][[t]]$ with the product as the group operation. Let $\S$ be either the space $\LL_0$ of arcs, or $\LL_0^*/\C^*$, or the space $\B_0$ of branches. To emphasize that the integration is with respect to the multiplicative structure, we shall denote such an integral by
$$
\int_\S (1-t^{v_f})^{-d\chi_g}.
$$
If $X_n=\{\varphi\in\S: v_f(\varphi)=n\}$,
then
$$
\int_\S t^{v_f}d\chi_g=\sum_{n=1}^\infty \chi_g(X_n)\cdot t^n.
$$
For the integral of $(1-t^{v_f})^{-1}$ one has
$$
\int_\S (1-t^{v_f})^{-d\chi_g}=\prod_{n=1}^\infty (1-t^n)^{-\chi_g(X_n)},
$$
where the expression $(1-t^n)^a$, $a\in K_0(\Var)[\LLL^{-1}]$, is understood in the sense of \cite{MRL}. There was constructed a, so called, power structure over the (semi)ring $\RR$ (any one of $S_0(\Var)$, $K_0(\Var)$, $K_0(\Var)[\LLL^{-1}])$.

\begin{definition}
A {\it power structure} over a (semi)ring $R$ is a map
$\left(1+t\cdot R[[t]]\right)\times {R} \to 1+t\cdot R[[t]]:$ $(A(t),m)\mapsto \left(A(t)\right)^{m}$, which possesses the properties:
\begin{enumerate}
\item $\left(A(t)\right)^0=1$,
\item $\left(A(t)\right)^1=A(t)$,
\item $\left(A(t)\cdot B(t)\right)^{m}=\left(A(t)\right)^{m}\cdot
\left(B(t)\right)^{m}$,
\item $\left(A(t)\right)^{m+n}=\left(A(t)\right)^{m}\cdot
\left(A(t)\right)^{n}$,
\item $\left(A(t)\right)^{mn}=\left(\left(A(t)\right)^{n}\right)^{m}$.
\end{enumerate}
\end{definition}

According to the construction in \cite{MRL}, for $a=-[Z]$, where $Z$ is a quasi-projective variety, one has
$$
(1-t)^a=\zeta_Z(t)=1+\sum_{k=1}^\infty [S^kZ]\cdot t^k,
$$
where $S^kZ=Z^k/S_k$ is the $k$th symmetric power of the space $Z$ ($\zeta_Z(t)$ is the Kapranov zeta function of the variety $Z$, see~\cite{L}), $(1-t)^{\frac{a}{\LLL^i}}=(1-t)^a\mbox{
\raisebox{-0.5ex}{$\vert$}}{}_{t\mapsto\frac{t}{\LLL^i}}
$\,.

In \cite{MRL} there also was defined a map $\Ex:t\cdot K_0(\Var)[\LLL^{-1}][[t]]\to
1+t\cdot K_0(\Var)[\LLL^{-1}][[t]]$ which was an isomorphism of the abelian groups (with addition and multiplication  as group operations respectively). It is defined by the equation
$$
\Ex\left(\sum\limits_{i=1}^\infty a_it^i\right)=\prod\limits_{i=1}^\infty(1-t^i)^{-a_i}.
$$
One can easily see that
\begin{equation}
\int_\S (1-t^{v_f})^{-d\chi_g}=\Ex\left(\int_\S t^{v_f}d\chi_g\right) \tag{*}.
\end{equation}

An intention to apply the constructions of this paper to the (multi variable) Alexander invariants of a collection of functions $f_1$, \dots, $f_r$ on $(\C^d,0)$ could lead to integrals of the form
$$
\int_\S (1-\tt^{\vv(\varphi)})^{d\chi_g}\,,
$$
where $\vv(\varphi)=(v_1(\varphi), \ldots, v_r(\varphi))$, $v_i(\varphi)=v_{f_i}(\varphi)$, $\tt=(t_1, \ldots, t_r)$, $\tt^\vv=t_1^{v_1}\cdot\ldots\cdot t_r^{v_r}$. Another reason to consider power series in several variables can be seen from Proposition~\ref{prop2} (if, say, $G=K_0(\Var)[[t]]$, $\psi(\varphi)=t^{v_f(\varphi)}$). This makes reasonable to consider expressions of the form $A(\tt)^M$, where
$A(\tt)=1+\sum\limits_{\kk\in\Z_{\ge0}^r,\,\kk\ne0}A_{\kk}\,\tt^{\kk}$,
$A_\kk\in\RR$, $M\in\RR$.
If there exists a power structure over a {\bf ring} $R$, there is a natural way to define the corresponding expression in the multi variable case as well. However, if $R$ is a {\bf semiring}, in general this does not work. Since elements of the Grothendieck semiring $S_0(\Var)$ of complex varieties have more geometric meaning (they are represented by "genuine" varieties, not by virtual ones), it is reasonable to give a geometric definition of this operation over this semiring.

It is possible (and convenient) to give the definition in a little bit more general setting. Let $S$ be an ordered abelian semigroup with zero such that each element $s\in S$ has only finitely many representations as sum of elements of $S$ (in particular, zero is the smallest element of $S$). For a (semi)ring $R$ there is defined the corresponding semigroup (semi)ring $R[[S]]$ which consists of formal sums (series) of the form $\sum\limits_{s\in S}r_s s$, where $r_s\in R$, with the natural operations:
$\sum r_s^{'} s + \sum r_s^{''} s=\sum (r_s^{'}+r_s^{''}) s$, $\left(\sum r_s^{'} s\right)\cdot \left(\sum r_s^{''} s\right)=\sum (r_{s_1}^{'}\cdot r_{s_2}^{''})(s_1+s_2)$, where in the last expression one should combine summands with the same $s_1+s_2$. Let $R_+[[S]]$ be the set (ideal) of series of the form $\sum\limits_{s\in S, s>0}r_s s$. The (semi)ring $R[[\tt]]$ of formal power series in $r$ variables $\tt=(t_1, \ldots, t_r)$ with coefficients from $R$ is the semigroup (semi)ring $R[[S]]$ for $S=\Z_{\ge0}^r$. 

It is convenient to describe the power structure over the Grothendieck semiring $S_0(\Var)$ in terms of graded spaces (sets). A graded space (with grading from $S_{>0}$) is a space $\A$ with a function $I_\A$ on it with values in $S_{>0}$. The element $I_\A(a)$ of the semigroup $S$ is called the weight of the point $a\in \A$. To a series $A\in 1+S_0(\Var)_+[[S]]$, $A=1+\sum\limits_{s\in S, s>0}[A_s]s$ one associates the graded space $\A=\coprod\limits_{s\in S, s>0} A_s$ with the weight function $I_\A$ which sends all points of $A_s$ to $s\in S$. In the other direction, to a graded space $(\A, I_\A)$ there corresponds the series $A=1+\sum\limits_{i=1}^\infty[A_s]s$ with $A_s=I_A^{-1}(s)$. To describe the series $A^{[M]}$, we shall describe the corresponding graded space $\A^M$ first. The space $\A^M$ consists of pairs $(K,\varphi)$, where $K$ is a finite subset of (the variety) $M$ and $\varphi$ is a map from $K$ to the graded space $\A$. The weight function $I_{\A^M}$ on $\A^M$ is defined by $I_{\A^M}(K,\varphi)=\sum\limits_{k\in K}I_\A(\varphi(k))$. This gives a set-theoretic description of the series $A^{[M]}$. To describe the coefficients of this series as elements of the Grothendieck semiring $S_0(\Var)$, one can write it as
$$
A^{[M]}=1+
\sum_{s_0\in S, s_0>0}
\left\{
\sum_{\kk:\sum sk_s=s_0}
\left[
\left(
(\prod_s M^{k_s})
\setminus\Delta
\right)
\times\prod_s A_s^{k_s}/\prod_s S_{k_s}
\right]
\right\}
\cdot s_0,
$$
where $\kk=\{k_s: s\in S, s>0, k_s\in\Z_{\ge0}\}$, $\Delta$ is the "large diagonal" in $M^{\Sigma k_s}$ which consists of $(\sum k_s)$-tuples of points of $M$ with at least two coinciding ones, the permutation group $S_{k_s}$ acts by permuting corresponding $k_s$ factors in
$\prod\limits_s M^{k_s}\supset (\prod_s M^{k_s})\setminus\Delta$ and the spaces $A_s$ simultaneously (the connection between this formula and the description above is clear).

\section{Motivic versions of the monodromy zeta function.}
Let 
$$
\eta_\S(t):=\int_\S(1-t^{v_f})^{-d\chi_g},
$$
where $\S=\LL_0$, $\LL_0^*/\C^*$, or the space of branches $\B_0$. Let us compute the specialization of the series $\eta_\S(t)$ under the (usual) Euler characteristic morphism.

\begin{proposition}\label{prop5}
For $\S=\LL_0$, $\chi(\eta_\S(t))=1$; for $\S=\LL_0^*/\C^*$ or $\B_0$,
$$
\eta(t)=\chi(\eta_\S(t))=\prod\limits_{k=1}^\infty\zeta_f(t^k),
$$
where $\zeta_f(t)$ is the classical monodromy zeta function of the germ $f$.
\end{proposition}

\begin{proof}
For $\S=\LL_0$, this follows from the expression for the integral
$\int_{\LL_0} t^{v_f} d\chi_g$ in terms of a resolution since all terms in it are divisible by $(\LLL-1)$. For $\S=\LL_0^*/\C^*$ or $\B_0$, from Propositions 3 and 4 it follows that under the Euler characteristic  morphism all terms corresponding to non trivial intersections of complements (i.e. to $I$ with $|I|>1$) vanish and
$$
\eta(t)=\chi(\eta_\S(t))=\prod_{i\in I_0^{'}} \prod_{k=1}^\infty (1-t^{k N_i})^{-\chi(\stackrel{\circ}{E_i})}
=\prod_{k=1}^\infty \left( \prod_{i\in I_0^{'}} (1-t^{k N_i})^{-\chi(\stackrel{\circ}{E_i})} \right)
$$
(this follows from the equation $\chi\left(A(t)^{[M]}\right)=\left(\chi(A(t))\right)^{\chi(M)}$, see \cite{MRL}). According to the A'Campo formula \cite{A'C} it is equal to
\begin{equation}
\eta(t)=\prod_{k=1}^\infty \zeta_f(t^k). \tag{**}
\end{equation}
\end{proof}

Let $\mu(n)$ be the M\"obius function:
$$
\mu(n)=
\begin{cases}
0&\mbox{ if $n$ has one or more repeated prime factors,}\\
1&\mbox{ if $n=1$,}\\
(-1)^k&\mbox{ if $n$ is a product of $k$ distinct primes.}
\end{cases}
$$

From the property  $\mu(1)=1$, $\sum\limits_{i\vert n}\mu(i)=0$ for $n>1$ one has the following statement.

\begin{corollary}
$$
\zeta_f(t)=\prod\limits_{i=1}^\infty(\eta(t^i))^{\mu(i)}.
$$
\end{corollary}

This leads to the following motivic versions of the monodromy zeta function.

\begin{definition}
For $\S=\LL_0^*/\C^*$ or $\B_0$, the series
$$
\zeta_{f,\S}(t):=\prod\limits_{i=1}^\infty(\eta_S(t^i))^{\mu(i)}
$$
will be called the $\mbox{\em arcs}/\C^*$ and {\em branches} motivic version of the monodromy zeta function for $\S=\LL_0^*/\C^*$ and $\S=\B_0$ respectively.
\end{definition}

In other words
\begin{eqnarray*}
\zeta_{f,\S}(t)&=&\Ex\left(\int_\S\sum\limits_{i=1}^\infty \mu(i)t^{iv_f}\,d\chi_g\right)=
\Ex\left(\sum\limits_{n=1}^\infty \left(\sum_{k|n} \mu(k) \chi_g(X_{n/k}) \right) t^{n}\,\right)\\
&=&\prod_{n\geq 1} (1-t^n)^{-\left(\sum_{k|n} \mu(k) \chi_g(X_{n/k}) \right)}\,.
\end{eqnarray*}

\begin{proposition}\label{prop6}
The series $\zeta_{f,\S}(t)\in K_0(\Var)[\LLL^{-1}][[t]]$ is an invariant of the germ $f$ and its specialization under the Euler characteristic morphism coincides with (the Taylor expansion of) the monodromy zeta function $\zeta_f(t)$.
\end{proposition}

Propositions \ref{prop3} and \ref{prop4} give:

\begin{theorem}\label{theo}
For a resolution $\pi:(\X,\D)\to(\C^d,0)$ of the germ $f$ one has:
$$
\zeta_{f,\LL_0^*/\C^*}(t)=\prod\limits_{m=1}^\infty\quad\prod\limits_{I\subset I_0,\, \emptyset\ne I\not\subset I''_0}\left(\prod\limits_{\{k_i\vert\,i\in I\}}(1-\LLL^{-\kk\,\nuunder}t^{m\kk\,\NN}) \right)^{-\mu(m)(\LLL-1)^{\#I-1}[\stackrel{\circ}{E_I}]}
$$
$($here $\kk=\{k_i\vert\,i\in I\}$, $\nuunder=\{\nu_i\vert\,i\in I\}$, $\NN=\{N_i\vert\,i\in I\}$, $\kk\,\nuunder=\sum\limits_{i\in I}k_i\nu_i$, \dots$)$.
\newline
If the resolution $\pi$ factorizes through the blowing-up $\pi_0:(\X_0,\C\P^{d-1})\to(\C^d,0)$ at the origin in $\C^d$, then:
$$
\zeta_{f,\B_0}(t)=\prod\limits_{m=1}^\infty\quad\prod\limits_{I\subset I_0,\, \emptyset\ne I\not\subset I''_0}\left(\prod\limits_{\{k_i\vert\,i\in I\}}(1-\LLL^{-\kk(\nuunder+\MM)}t^{m\kk\,\NN}) \right)^{-\mu(m)(\LLL-1)^{\#I-1}[\stackrel{\circ}{E_I}]}
$$
$($$\MM=\{M_i\vert\,i\in I\}$$)$.
\end{theorem}

\section{Final remarks.}
The function $\int_{\LL_0} t^{v_f}\,d\chi_g$ is the "naive" zeta function of J.~Denef and F.~Loeser (see \cite{DL3}) which is a rational function (this follows from its description in terms of a resolution). One of the most interesting problems about this function is the "Monodromy Conjecture" which states that there is a set $S= \{(\nu,N):\, \nu,N\in\N,N>0 \}$ such that the "naive" zeta function always belongs to $K_0(\Var)[\LLL^{-1}][(1-\LLL^{-\nu}t^N)^{-1}]_{\{(\nu,N)\}\in S}$
and if $q=-\nu/N,(\nu,N) \in S$, then $\exp(-2 i\pi q)$ is an eigenvalue of the classical local monodromy operator around zero at some point $P\in f^{-1}(0)$. 

After Proposition~\ref{prop2} it is clear that the Monodromy Conjecture for $
\int_{\LL_0} t^{v_f}\,d\chi_g$ is equivalent the Monodromy Conjecture for $
\int_{\LL_0^*/\C^*} t^{v_f}\,d\chi_g$. Because of the identities $(*)$ and $(**)$ one has   
\begin{eqnarray*}
\zeta_{f,\S}(t)=\Ex\left(\sum\limits_{i=1}^\infty \mu(i) \int_{\LL_0^*/\C^*} t^{iv_f}\,d\chi_g\right), \\
\chi\left(\Ex\left(\int_{\LL_0^*/\C^*} t^{v_f}d\chi_g\right)\right)=\prod_{k=1}^\infty \zeta_f(t^k).
\end{eqnarray*}
These seems to be the first direct formulae which connects the Denef--Loeser zeta function and the classical monodromy zeta function.

The Monodromy Conjecture was originally stated in the $p$-adic
case for the Igusa local zeta functions, see e.g. \cite{de:91}. In \cite{dl:92}
J.~Denef and F.~Loeser introduced an analytic invariant called \emph{local topological zeta function} of a germ $f$  (as a kind of limit of the
local Igusa zeta function) whose initial definition was written in terms of a resolution:
$$
Z_{top,0}(f,s)= \sum_{I\subset I_0, \emptyset\ne I\not\subset I_0^{''}}
\chi(\stackrel{\circ}{E_I}) \prod_{i\in I} \frac{1}{N_is+\nu_i}
$$
besides that it does not depend on it. If first we substitute $t$ by $\LLL^{-s}$ in $\int_{\LL_0} t^{v_f}\,d\chi_g$,
then expand $\LLL^{-s}$ and $(\LLL-1)(1-\LLL^{-\nu+N s})^{-1}$
into series in $\LLL-1;$ and finally take the Euler characteristic
then one gets $Z_{top,0}(f,s)$. Also the Monodromy Conjecture has been stated for this function in \cite{dl:92}. See \cite{ve} for more information about the conjecture.

In the case of integration over branches one can do the same procedure to get a new analytic invariant of the germ $f$. The function 
$$
Z_{\B,0}(f,s):=\chi((\LLL-1)\int_{\B}t^{v_f}\,d\chi_g
\raisebox{-0.5ex}{$\vert_{{T\mapsto\,\LLL^{-s}}}$})
$$
is rational and in terms of a resolution which factorizes through the blowing-up $\pi_0:(\X_0,\C\P^{d-1})\to(\C^d,0)$ at the origin in $\C^d$ one has
$$
Z_{\B,0}(f,s)= \sum_{I\subset I_0, \emptyset\ne I\not\subset I_0^{''}}
\chi(\stackrel{\circ}{E_I}) \prod_{i\in I} \frac{1}{N_is+\nu_i+M_i}\,.
$$
In particular one could ask if the Monodromy Conjecture holds for $Z_{\B,0}(f,s)$. The example $f(x,y)=(x^2+y^3)(y^2+x^3)$ shows that this is not the case. 
The monodromy zeta function $\zeta_f(t)$ is $(1-t^{10})^2/(1-t^5)^2$ while 
$$ 
Z_{\B,0}(f,s)=\frac{8s^2+24s+14}{(10s+7)(4s+3)(1+s)}\,.
$$

\end{document}